\documentclass[12pt]{article}
\usepackage{amsmath,amssymb,amsfonts,amsthm}
\usepackage[all]{xy}

\newcounter{item}[section]
\newcounter{kirshr}
\newcounter{kirsha}
\newcounter{kirshb}
\newenvironment{enumroman}{\setcounter{kirshr}{1}
\begin{list}{(\roman{kirshr})}{\usecounter{kirshr}} }{\end{list}}
\newenvironment{enumarab}{\setcounter{kirshb}{1}
\begin{list}{(\arabic{kirshb})}{\usecounter{kirshb}} }{\end{list}}
\newenvironment{athm}[1]{\vskip3mm\par\noindent
{\bf #1 }. \slshape }
{\upshape\par\vskip10pt minus3pt}
\newtheorem{theorem}{Theorem}[section]

\newtheorem{lemma}[theorem]{Lemma}

\newenvironment{demo}[1]{\noindent{\bf #1.}\upshape\mdseries}
{\nopagebreak{\hfill\rule{2mm}{2mm}\nopagebreak}\par\normalfont}
\theoremstyle{definition}
\newtheorem{remark}[theorem]{Remark}

\newtheorem{definition}[theorem]{Definition}

\def\C{{\mathfrak{C}}}

\def\Nr{{\mathfrak{Nr}}}
\def\Fr{{\mathfrak{Fr}}}
\def\Sg{{\mathfrak{Sg}}}

\def\A{{\mathfrak{A}}}
\def\B{{\mathfrak{B}}}
\def\C{{\mathfrak{C}}}
\def\D{{\mathfrak{D}}}

\def\CA{{\bf CA}}

\def\K{{\bf K}}
\def\K{{\bf K}}

\def\RCA{{\bf RCA}}

\def\Rd{{\mathfrak{Rd}}}
\def\(R)RA{{\bf (R)RA}}

\def\R{\mathbb{R}}

 \def\CA{{\sf CA}}
\def\B{{\sf B}}

\def\K{{\sf K}}
 
\def\Nr{{\mathfrak{Nr}}}

\def\Nr{{\mathfrak{Nr}}}

\def\A{{\mathfrak{A}}}
\def\B{{\mathfrak{B}}}
\def\C{{\mathfrak{C}}}
\def\D{{\mathfrak{D}}}

\def\CA{{\bf CA}}

\def\RCA{{\bf RCA}}
\def\RQEA{{\bf RQEA}}

\def\R{\cal{R}}

\def\Nr{{\mathfrak{Nr}}}
\def\Fr{{\mathfrak{Fr}}}
\def\Sg{{\mathfrak{Sg}}}

\def\Rd{{\mathfrak{Rd}}}

\def\CA{{\bf CA}}
\def\RCA{{\bf RCA}}
\def\K{{\bf K}}

\def\(R)RA{{\bf (R)RA}}

\def\R{\mathbb{R}}

\def\Nr{{\mathfrak{Nr}}}
\def\Fr{{\mathfrak{Fr}}}
\def\Sg{{\mathfrak{Sg}}}

\def\Rd{{\mathfrak{Rd}}}

\def\CA{{\bf CA}}
\def\RCA{{\bf RCA}}
\def\K{{\bf K}}

\def\(R)RA{{\bf (R)RA}}

\def\R{\mathbb{R}}

\def\R{\mathbb{R}}

\def\R{\mathbb{R}}

 \def\CA{{\sf CA}}

\def\K{{\bf K}}

\def\A{{\mathfrak{A}}}
\def\B{{\mathfrak{B}}}
\def\C{{\mathfrak{C}}}
\def\D{{\mathfrak{D}}}
\def\P{{\mathfrak{P}}}

\def\Nr{{\mathfrak{Nr}}}
\def\F{{\mathfrak{F}}}
\def\CA{{\bf CA}}
\def\RCA{{\bf RCA}}

\def\Sg{{\mathfrak Sg}}

\title{Representability, amalgamation in connection to the finitizability problem for  Heyting polyadic algebras\\ 
The Finitizability problem for predicate intuitionistic logic}
\author{Tarek Sayed Ahmed}

\begin{document}
\maketitle

\begin{abstract}
We prove completeness and interpolation for infinitary predicates intuitionistic logic, 
by proving new neat embedding theorems for Heyting polyadic algebras.
Also we formulate and solve the analogue of the lengthy discussed
finitizability problem in (classical) algebraic logic, which has to do with finding finitely based varieties adequate for the algebraisation
of predicate intuinistic logic circumventing non-finite axiomatizability results.

\footnote{ 2000 {\it Mathematics Subject Classification.} Primary 03G15.

{\it Key words}: algebraic logic, neat reducts, cylindric algebras, amalgamation} 


\end{abstract}

We follow standard notation. $G$ will always be a semigroup of transformations
on a set or ordinal (well ordered set) under the operation of composition.
All undefined terms, concepts in this part of the paper which we use extensively without warning  
are found in Part 1.  Cross references will be given.

This part is organized as follows. 
In section 2  we give a simple proofs  of representation theorems for the countable case, using ideas of Ono. 
In section 2, we prove various representation 
theorems in the presence of diagonal elements.

\section{Simpler proofs}

Throughout $G$ is a semigroup. We start by giving a simpler proof of our representability results when $G$ consists only of finite transformations
and $G$ is strongly rich. We assume that everything is countable; algebras, dimensions and semigroups.
We do not study the case when $G$ consists of all transformations since in this case we have uncountably many operations
(even if the dimension is countable. 

Dimensions will be specified by countable 
ordinals for some time to come. An {\it $\omega$ dilation} of an algebra of dimension $\alpha$ is its minimal dilation of  dimension $\alpha+\omega.$
Our next  proof is a typical Henkin construction that does not involve 
zigzagging. It as an algebraic version of a proof of Ono, but it proves much more.

We treat the two cases simultaneously; $V$ denotes the class of algebras and 
$\Fr_{\beta}^{\rho}V$ denotes the free countable $V$ algebra,
where $\rho$ is not dimension restricting in the second  case while $\alpha\sim \rho(i)$ is infinite in 
the first  case [cf. definition ? in part 1]. 

\begin{definition}
\begin{enumarab}
\item Let $\A$ be an algebra. Let $\Gamma, \Delta\subseteq \A$. We write $\Gamma\to \Delta$, 
if there exist $n,m\in \omega$, $a_1,\ldots a_n\in \Gamma$, $b_1,\ldots b_m\in \Delta$ such that
$a_1\land \ldots a_n\leq b_1\lor\ldots b_m$. Recall that a  pair $(\Gamma, \Delta)$ is consistent if not $\Gamma\to \Delta.$

\item We say that a theory $(\Gamma, \Delta)$ is complete, if it is consistent and $\Gamma\cup \Delta=\A.$
\item A theory $(\Gamma, \Delta)$ is Henkin complete in $\A$ if it is complete and saturated. 
$\Gamma\subseteq A$ is saturated, if there exists (equivalently for all) $\Delta\subseteq A$,
such that $(\Gamma, \Delta)$ is saturated.
\end{enumarab}
\end{definition}

\begin{definition}
\begin{enumarab}
\item Let $\Gamma\subseteq \Sg^{\A}X_1$ and $\Theta,\Lambda\subseteq \Sg^{\A}X_2$. Then $a\in \Sg^{\A}(X_1\cap X_2)$ separates $\Gamma$ from
$(\Theta, \Delta)$ if $\Gamma\to \{a\}$ and $\Theta\cup \{a\}\to \Lambda$.  In this case we say that $\Gamma$ can be separated from 
$(\Theta, \Delta)$. Otherwise $\Gamma$ is inseparable
from $(\Theta, \Delta)$ with respect to $\A$, $X_1$, $X_2$,  or simply inseparable, when $\A$, $X_1$ and $X_2$ are clear from context.

\item We say that $a\in \Sg^{\A}(X_1\cap X_2)$ separates $\Gamma$ from
$\Delta$ if $a$ separates $\Gamma$ from $(\emptyset, \Delta),$ that is to say, if $\Gamma\to \{a\}$ and $\{a\} \to \Delta$.  
\end{enumarab}
\end{definition}

\begin{lemma}\label{Ono}(Essentially Ono's) let $G$ be the semigroup of all finite transformations on $\alpha$. Let $\A\in GPHA_{\alpha}$
and assume further that $\alpha\sim \Delta x$ is infinite for every $x$ in $A$. Suppose that $\Gamma\subseteq \Sg^{\A}X_1$ and 
$\Theta, \Lambda\subseteq \Sg^{\A_2}X_2$. If $\Gamma$ is inseparable from $(\Theta,\Lambda)$
with respect to $X_1$ and $X_2$, then there exist an $\omega$ dilation $\B$ of $\A$, 
$\Gamma'\subseteq \Sg^{\B}X_1$ 
and $\Theta', \Delta'\subseteq \Sg^{\B}X_2,$ such that
\begin{enumerate}
\item $\Gamma\subseteq \Gamma'$ and $\Gamma'$ is saturated in $\Sg^{\B}X_1$.
\item $\Theta\subseteq  \Theta'$ and $\Lambda\subseteq \Lambda'$ and $(\Theta', \Lambda')$ is Henkin complete.
\item $(\Gamma'\cap \Theta', \Delta')$ is Henkin complete and $\Delta'=\Lambda'\cap \Sg^{\B}(X_1\cap X_2)$.
\item $\Gamma'$ is inseparable from $(\Theta', \Lambda').$
\end{enumerate}
\end{lemma}
\begin{demo}{Proof} For each $i\in \omega$, we will construct subsets $\Gamma _{i}$
of $\Sg^{\B}(X_1)$ and $\Theta _{i}$, $\Lambda _{i}$ of $\Sg^{\B}(X_2)$ satisfying the following 

\begin{description}
\item  (1) \qquad $\Gamma =\Gamma _{0}\subseteq \Gamma _{1}\subseteq \Gamma
_{2}\subseteq ...$,

\item   \qquad $\quad \Theta =\Theta _{0}\subseteq \Theta _{1}\subseteq
\Theta _{2}\subseteq ...,$

\item   \qquad \quad $\Lambda =\Lambda _{0}\subseteq \Lambda _{1}\subseteq
\Lambda _{2}\subseteq ...$.


\item  (2) \qquad $\Gamma _{i}$ is inseparable from $(\Theta _{i},\Lambda
_{i})$ with respect to $\Sg^{\B}(X_1\cap X_2)$.
\end{description}
Let $(a_i: i\in \omega)$ be an enumeration of  $\Sg^{\B}(X_1\cap X_2)$, $(b_i: i\in \omega)$ 
be an enumeraton of $\Sg^{\B}(X_1)\sim \Sg^{\B}(X_1\cap X_2)$ and
$(c_i: i\in \omega)$ be an enumeration of $\Sg^{\B}(X_2)\sim \Sg^{\B}(X_1\cap X_2)$.

First observe  that $\Gamma $ is inseparable from $(\Theta ,\Lambda )$
with respect to the big algebra $\B$ and $X_1$ and $X_2$.  For if not, then there exists $b\in \Sg^{\B}(X_1\cap X_2)$ 
such that $\Gamma\to \{b\}$ and $\Theta\cup \{b\}\to \Lambda$. By lemma \ref{cylindrify}, there is a finite set $J,$ such that 
$${\sf c}_{(J)}b\in \Nr_{\alpha}\Sg^{\B}(X_1\cap X_2)=\Sg^{\Nr_{\alpha}\B}(X_1\cap X_2)=\Sg^{\A}(X_1\cap X_2),$$ 
but then we would have 
$\Gamma\to \{{\sf c}_{(J)}b\}$ and $\Theta\cup \{ {\sf c}_{(J)}b\}\to \Lambda$ which is impossible. So,  (2)  holds when $%
i=0$. Next, suppose that $i>0$ and inductively for each $k<i$, $\Gamma_{k},\Theta _{k}$ and $\Lambda_{k}$ 
satisfying the conditions (1), (2), are defined.
We also assume inductively that
$$\beta\sim \bigcup_{x\in \Gamma_i}\Delta x \cup\bigcup_{x\in \Theta_i}\Delta x\cup\bigcup_{x\in \Lambda_i}\Delta x$$
is infinite for all $i$. This is clearly satisfied for the base of induction.
We treat three cases separately.
\begin{enumerate}
\item $i=3m+1$. Since the sentence $a_{m}$ belongs to $\Sg^{\B}(X_1\cap X_2)$
either $\Gamma _{i-1}\cup \{a_{m}\}$ is inseparable from $%
(\Theta _{i-1}\cup \{a_{m},\Lambda _{i-1})$, or $\Gamma _{i-1}$ is
inseparable from $(\Theta _{i-1},\Lambda _{i-1}\cup \{a_{m}\})$.
Suppose that the former case holds. 

We distinguish between two subcases:

{\bf Subcase (1)}: $a_m\neq {\sf c}_kx$ for all $k<\beta$ and $x\in B$.
Set $\Gamma _{i}=\Gamma _{i-1}\cup
\{a_{m}\}$, $\Theta _{i}=\Theta _{i-1}\cup \{a_{m}\}$ and $\Lambda
_{i}=\Lambda _{i-1}$. 

{\bf Subcase (2)}: $a_m={\sf c}_ka$ for some $k<\beta$ and $a\in B$.
Then set $\Gamma _{i}=\Gamma _{i-1}\cup \{a_{m}, {\sf s}_j^ka\}$, $\Theta
_{i}=\Theta _{i-1}\cup \{a_{m}, {\sf s}_j^ka\}$ and $\Lambda
_{i}=\Lambda _{i-1}$, where
$$j\in \beta\sim \bigcup_{x\in \Gamma_{i-1}}\Delta x\cup\bigcup_{x\in \Theta_{i-1}}\Delta x\cup\bigcup_{x\in \Lambda_{i-1}}\Delta x.$$ 
Such a $j$ exists by the induction hypothesis. Then $\Gamma _{i}$ is inseparable
from $(\Theta _{i},\Lambda _{i})$. 

Suppose next that the latter case holds. Then, define $\Gamma _{i}=\Gamma _{i-1}$, $%
\Theta _{i}=\Theta _{i-1}$ and $\Lambda _{i}=\Lambda _{i-1}\cup \{a_{m}\}$.

In each subcase, as easily checked, the conditions (1) and (2) hold for $\Gamma _{i}$, $\Theta
_{i}$ and $\Lambda _{i}$.

\item The case where $i=3m+2$. Suppose first that some sentence in $\Sg^{\B}(X_1\cap X_2)$ that 
separates $\Gamma _{i-1}\cup \{b_{m}\}$ from $(\Theta
_{i-1},\Lambda _{i-1})$. Then, define $\Gamma _{i}=\Gamma _{i-1}$, $\Theta
_{i}=\Theta _{i-1}$, $\Lambda _{i}=\Lambda _{i-1}$. Next suppose otherwise.

{\bf Subcase 1}:
When $b_{m}$ is not of the form ${\sf c}_kb$ for any $k\in \beta$ and $b\in B$; we define $\Gamma
_{i}=\Gamma _{i-1}\cup \{b_{m}\}$, $\Theta _{i}=\Theta _{i-1}$ and $\Lambda
_{i}=\Lambda _{i-1}$. 

{\bf Subcase 2}: When $b_{m}$ is equal to ${\sf c}_kb$ for some $k\in \beta$  and $b\in B$,  define $\Gamma _{i}=\Gamma _{i-1}\cup \{b_m,  {\sf s}_j^kb\}$, $%
\Theta _{i}=\Theta _{i-1}$ and $\Lambda _{i}=\Lambda _{i-1}$ where  
$$j\in \beta\sim \bigcup_{x\in \Gamma_{i-1}}\Delta x\cup\bigcup_{x\in \Theta_{i-1}}\Delta x\cup\bigcup_{x\in \Lambda_{i-1}}\Delta x.$$ 
In any case, we have $\Gamma_i$ is inseparable from $(\Theta _{i},\Lambda _{i}).$

\item The case where $i=3m+3$. As in the two previous cases, either $\Gamma _{i-1}$ is
inseparable from $(\Theta _{i-1}\cup \{c_{m}\}, \Lambda_{i-1})$ or $\Gamma _{i-1}$
is inseparable from $(\Theta _{i-1},\Lambda _{i-1}\cup \{c_{m}\})$. If it is
the latter case, then set $\Gamma _{i}=\Gamma _{i-1}$, $\Theta _{i}=\Theta
_{i-1}$ and $\Lambda _{i}=\Lambda _{i-1}\cup \{c_{m}\}$. 
If it is the former case,  when $c_{m}$ is not of the form ${\sf c}_kc$ for any $k\in \beta$ and $c\in B$,  define $\Gamma
_{i}=\Gamma _{i-1}$, $\Theta _{i}=\Theta _{i-1}\cup \{c_{m}\}$ and $\Lambda
_{i}=\Lambda _{i-1}$. If there exists $c\in \B$, and $k\in \beta$ such that $c_{m}= {\sf c}_kc,$ set $\Gamma _{i}=\Gamma _{i-1}$, 
$\Theta _{i}=\Theta _{i-1}\cup
\{c_m , {\sf s}_j^kc\}$ and $\Lambda _{i}=\Lambda _{i-1}$ where
$$j\in \beta\sim \bigcup_{x\in \Gamma_{i-1}}\Delta x\cup\bigcup_{x\in \Theta_{i-1}}\Delta x\cup\bigcup_{x\in \Lambda_{i-1}}\Delta x.$$ 

In each case, the conditions  (1) and (2)  are satisfied.

\end{enumerate}
Now, we define $\Gamma ^{\prime }=\bigcup_{i<\omega }\Gamma _{i}$, $\Theta
^{\prime }=\bigcup_{i<\omega }\Theta _{i}$ and $\Lambda ^{\prime}=\bigcup_{i<\omega }\Lambda _{i}$. Clearly, $\Gamma ^{\prime }$ is
inseparable from $(\Theta ^{\prime },\Lambda ^{\prime })$, by (1) and (2).

We will prove that $\Gamma ^{\prime }$ is saturated in $\Sg^{\B}(X_1)$.

We first prove

{\bf Claim} For any  
$b$ in $\Sg^{\B}X_1$, $b$ is in $\Gamma ^{\prime }$ if and
only if $\Gamma _{i}\cup \{b\}$ is inseparable from $(\Theta _{i},\Lambda
_{i})$ for every $i$. 
 
\begin{demo}{Proof of Claim} Suppose that for each $i,$ $\Gamma _{i}\cup \{b\}$
is inseparable from $(\Theta _{i},\Lambda _{i})$. Assume that  $b=b_{m}$ 
and let $k=3m+1$ if $b$ belongs to $\Sg^{\B}(X_1\cap X_2)$  and $k=3m+2$ otherwise. 
Then $b\in \Gamma _{k}\subseteq \Gamma^{\prime }$, by the fact that $\Gamma _{k-1}\cup \{b_{m}\}$ is inseparable
from $(\Theta _{k-1},\Lambda _{k-1})$. 

Conversely, suppose that $b\in \Gamma ^{\prime }$. Let $J=\{i\in \omega: b\in \Gamma_i\}$ and let $h= min J.$ 
Then for any $i$ such
that $h\leq i,\Gamma _{i}\cup \{b\}$ is inseparable from $(\Theta_{i},\Lambda _{i})$, since $\Gamma _{i}\cup \{b\}=\Gamma _{i}$. 
So, suppose, seeking a contradiction  that for some $i<h$,  we have $\Gamma _{i}\cup \{b\}$ is separated from $(\Theta_{i},\Lambda _{i})$. 
Then, it follows that $\Gamma _{h}(=\Gamma _{h}\cup \{b\})$ is also separated from $(\Theta _{h},\Lambda _{h})$. But,
this contradicts (2). Hence, $\Gamma _{i}\cup \{b\}$ is inseparable from 
$(\Theta _{i},\Lambda _{i})$ for each $i$.
\end{demo}
Now, we will show that $(\Gamma ^{\prime },\Sg^{\B}X_1\sim \Gamma^{\prime })$ is consistent. If not, then, there exist $%
a_{1}\ldots a_{m}$ in $\Gamma ^{\prime }$ and $b_{1}\ldots b_{n}$ in $\Sg^{\B}(X_1)\sim \Gamma$ such that 
$\bigwedge_{j=1}^{m}a_j\leq  \bigvee_{j=1}^{n}b_{j}$. For each $j,$ there exist a
number $k_{j}$ and an element $d_{j}\in \Sg^{\B}(X_1\cap X_2),$  such
that $d_{j}$ separates $\Gamma _{k_j}\cup \{b_{j}\}$ from $(\Theta _{k_j},\Lambda _{k_j})$. Take $k\in \omega$,  such that $k_{j}<k$ for each 
$j$ and $a_{1}\ldots a_{m}\in \Gamma _{k}$. Such a $k$ of course exists.  Now $\bigvee_{j=1}^{n}d_{j}$ in $\Sg^{\B}(X_1\cap X_2)$, 
separates $\Gamma_{k}\cup \{\bigvee_{j=1}^{n}b_{j}\}$ from $(\Theta _{k},\Lambda _{k})$.
Since $a_{1}\ldots a_{m}\in \Gamma _{k}$ and $\bigwedge_{j=1}^{m}a_{j}\leq
\bigvee_{j=1}^{n}b_{j}$, then we can infer that $\bigvee_{j=1}^{n}d_{j}$ separates
also $\Gamma _{k}$ from $(\Theta _{k},\Lambda _{k})$. But this contradicts
(2). 

Hence, $(\Gamma ^{\prime }, \Sg^{\B}X_1\sim \Gamma ^{\prime })$
is consistent. By construction it is also
Henkin-complete in $Sg^{\B}X_1$. Thus, $\Gamma ^{\prime }$ is
saturated in $\Sg^{\B}(X_1)$.
That $(\Theta ^{\prime },\Lambda ^{\prime })$ is
Henkin-complete in $\Sg^{\B}X_2$ is absolutely straightforward. From this, it immediately
follows that $(\Gamma ^{\prime }\cap \Theta ^{\prime },\Delta ^{\prime })$
is consistent. Let $a\in \Sg^{\B}(X_0\cap X_1)$. Suppose
that $a\notin \Delta ^{\prime }$. Then, $a\in \Theta ^{\prime }
$, since $(\Theta ^{\prime },\Lambda ^{\prime })$ is complete 
in $\Sg^{\B}X_1$. On the other hand, $a$ is in either  $\Gamma ^{\prime }
$ or  $\Lambda ^{\prime }$ by our construction. Hence, $a\in \Gamma
^{\prime }$ and, therefore, $a\in \Gamma ^{\prime }\cap \Theta ^{\prime }$.
This means that $(\Gamma ^{\prime }\cap \Theta ^{\prime },\Delta ^{\prime })$
is complete, and also Henkin-complete by construction.
\end{demo}
\begin{remark} The proof works verabtim when $G$ is stronly rich.
\end{remark}.
\begin{theorem}\label{Ono2} Let $\A=\Fr_{\mu}^{\rho}V$, where $|\mu|=\omega$, $\Gamma_0\subseteq \Sg^{\A}X_1$ 
and $\Theta_0, \Lambda_0\subseteq \Sg^{\A}X_2$. If $\Gamma_0$ is inseperable from $(\Theta_0, \Lambda_0)$ 
then the theory  $(\Gamma_0\cup \Theta_0, \Lambda_0)$ is satisfiable. That is to say, 
there exists $\mathfrak{K}=(K,\leq \{X_k\}_{k\in K}\{V_k\}_{k\in K}),$ a homomorphism $\psi:\A\to \mathfrak{F}_K,$
$k_0\in K$, and $x\in V_{k_0}$,  such that for all $p\in \Gamma_0\cup \Theta_0$ 
if $\psi(p)=(f_k)$, then $f_{k_0}(x)=1$ and for all $p\in \Lambda_0$
if $\psi(p)=(g_k)$, then $g_{k_0}(x)=0.$

\end{theorem}

\begin{demo}{Proof} We provide the proof when $G$ is the set of all finite transformations. Fix $\rho$ such that 
$\alpha\sim \rho(i)$ is infinite
for every $i\in \mu$. Form a sequence of minimal dilations 
$\A=\B_0\subseteq \B_1\subseteq\ldots \B_{\infty},$ where for  $n\leq \omega$, $\B_n\in GPHA_{\alpha+n}$.
Now by lemma \ref{dl}, we have for $k<l\leq \omega$, $\B_k=\Nr_k\B_l$ and for all $X\subseteq B_k$, we have 
$\Sg^{\B_k}X=\Nr_k\Sg^{\B_l}X$.  In this case we can give a more concrete description of the dilations.
It will turn out, as we proceed to show in a minute, that $\B_l$  is the dimension restricted free algebra on $\mu$ generators
in $\alpha+l$ dimensions restricted by $\rho$. For $n\leq \omega$, let $\beta=\alpha+n$. 
For brevity, we write $Hom(\A, \B)$ for the set of all homomorphisms from $\A$ to $\B$. Let $V_I=GPHA_I$.
In view of  lemma \ref{dl}, it suffices to show that 
the sequence $\langle \eta/Cr_{\mu}^{\rho}V_{\beta}:\eta<\mu\rangle$
$V_{\alpha}$ - freely generates
$\Nr_{\alpha}\Fr_{\mu}^{\rho}(V_{\beta})$.  
Towards this end, let $\B\in
V_{\alpha}$ and $a=\langle a_{\eta} :\eta<\mu\rangle\in {}^{\mu}B$
be such that $\Delta a_{\eta}\subseteq \rho(\eta)$ for all
$\eta<\mu$.  Since $\Fr_{\mu}^{\rho}(V_{\alpha})$ is in
$Dc_{\alpha}$, we have $|\alpha\setminus \bigcup_{\xi\in
\Gamma}\rho\xi|\geq \omega$ for each finite $\Gamma\subseteq \mu$.
Assuming that $Rg a$ generates $\B$, we have $\B\in Dc_{\alpha}$.
Therefore $\B\in S\Nr_{\alpha}V_{\beta}.$ Let
$\D=\Fr_{\mu}^{\rho}(V_{\beta})$. 

We claim that  $x=\langle
\eta/Cr_{\mu}^{\rho}V_{\beta}:\eta<\mu\rangle\in {}^{\mu}D$
$S\Nr_{\alpha}K_{\beta}$ - freely generates
$\Sg^{\Nr_{\alpha}\D}Rgx$.  Indeed, consider $\C\in
\Nr_{\alpha}V_{\beta}$ and $y\in {}^{\mu}C$ such that $\Delta
y_{\eta}\subseteq \rho \eta$ for all $\eta<\mu$. Let $\C'\in
V_{\beta}$ be such that $\C=\Nr_{\alpha}\C'$. Then clearly $y\in
{}^{\mu}C'$ and $\Delta y_{\eta}\subseteq \alpha$ for all
$\eta<\mu$. Then there exists $h\in Hom(\D, \C')$ such that $h\circ
x=y$. Hence $h\in Hom(\Rd_{\alpha}\D, \Rd_{\alpha}\C')$, thus $h\in
Hom(\Sg^{\Rd_{\alpha}\D}Rg x, \Sg^{\Rd_{\alpha}\C'}h(Rgx))$. Since
$Rg x\subseteq Nr_{\alpha} D$, we have $h\in
Hom(\Sg^{\Nr_{\alpha}\D}Rgx, \C)$. We have proved our claim.

Therefore there exists, by universal property of dimension restricted free algebras, a homomorphism 
$h:\Sg^{\Nr_{\alpha}\D}Rgx\to \B$ such that
$h(\eta/Cr_{\mu}^{\rho}K_{\beta})= a_{\eta}$. But
$\Nr_{\alpha}\Fr_{\mu}^{\rho}(K_{\beta})=\Nr_{\alpha}(\Sg^{\D}Rgx)=\Sg^{\Nr_{\alpha}\D}Rgx.$
Therefore  $\langle \eta/Cr_{\mu}^{\rho}(V_{\beta}):
\eta<\mu\rangle$ $V_{\alpha}$ freely generates
$\Nr_{\alpha}\Fr_{\mu}^{\rho}(V_{\beta})$.
We have proved that $\Fr_{\mu}^{\rho}V_{\alpha}\cong \Nr_{\alpha}\Fr_{\mu}^{\rho}V_{\beta}$, 
via an isomorhism that fixes the generators. 
Hence the (minimal) dilations are nothing more than dimension restricted free algebras in extra dimensions $\leq \omega$, also restriced
by $\rho$

Now  $\Gamma_0$ is inseparable from $(\Theta_0, \Lambda_0)$ with respect to $\B, X_1$ and $X_2$,
and so there exist by lemma \ref{Ono}, $\Gamma_0'\subseteq \Sg^{\B_1}(X_1)$, and $\Theta_0', \Lambda_0'\subseteq \Sg^{\B_1}(X_2)$ such that
$(\Gamma_0', \Sigma_0'), (\Theta_0', \Lambda_0')$ and $(\Gamma_0'\cap \Theta_0', \Delta_0')$ are Henkin complete 
in $\Sg^{\B_1}(X_1)$, $\Sg^{\B_2}(X_2)$ and $\Sg^{\B}(X_1\cap X_2)$ respectively,
where $\Sigma_0'=\Sg^{\B_1}(X_1)\sim \Gamma_0'$ and $\Delta_0'=\Lambda_0'\cap \Sg^{\B_1}(X_1\cap X_2).$

We define inductively for each $k\in \omega$ a set $M_k$ each of whose members is a pair of the form $(\Gamma, \theta)$ 
such that $\Gamma\subseteq \Sg^{\B_k}(X_1)$, and $\Theta\subseteq \Sg^{\B_k}(X_2),$ and the following hold:

\begin{enumarab}

\item $\Gamma_0'\subseteq \Gamma, \Theta_0'\subseteq \Theta$,

\item $(\Gamma, \Sigma), (\Theta, \Lambda)$ and $(\Gamma\cap \Theta, \Sigma\cap \Lambda)$ 
are Henkin complete in $\Sg^{\B_k}X_1$, $\Sg^{\B_k}(X_2)$ and $\Sg^{\B_k}(X_1\cap X_2)$,
respectively, where
$\Sigma=\Sg^{\B_k}(X_1)\sim \Gamma$ and $\Lambda=\Sg^{\B_k}(X_k)\sim \Theta$.

\end{enumarab}

1. At $k=0$, set $W_0=\{(\Gamma_0', \Theta_0')\}.$

2. Suppose that $W_{k-1}$ has already been defined and that each member of $W_{k-1}$ satisfies (1) and (2).
For each $(\Gamma,\Theta)\in M_{k-1}$ for all $a\in \Sigma\cup \Lambda$ of the form $d\to e$ or ${\sf q}_ia$ we construct a pair satisfying (1) and (2) 
and put in into $M_k$. Here (a) and (b) are analogous to Claims $1$ and $2$ in the proof of \ref{main}.

First take $d\to e\in \Lambda$. Then $\Gamma$ is inseperable from $(\Theta\cup \{d\},e)$ with respect to $\Sg^{B_{k-1}}(X_1\cap X_2)$, 
so there exists $\Gamma'\subseteq \Sg^{\B_k}(X_1)$ and $\Theta',\Lambda'\subseteq \Sg^{\B_k}(X_2)$ 
satisfying (1) and (2), and $\Theta\cup \{d\}\subseteq \Theta'$ and $e\in \Lambda'.$
If  ${\sf q}_ia\in\Lambda$. Then $\Gamma$ is inseprable from $(\Theta, {\sf s}_k^ia)$ where 
$k\in dimB_k\sim dimB_{k-1},$ with respect to $\Sg^{B_{k-1}}(X_1\cap X_2)$.
Now let
$K=\bigcup_{k\in \omega}W_k$; $K$ is the set of worlds.
If $i=(\Delta, \Gamma)\in W_k,$ let $M_k=\omega+k.n.$ 
If $i_1=(\Delta_1, \Gamma_1)$ and $i_2=(\Delta_2, \Gamma_2)$ are in $K$, then set
$$i_1\leq i_2\Longleftrightarrow  M_{i_1}\subseteq M_{i_2}, \Delta_1\subseteq \Delta_2. $$
This is a preorder on $K$; define the kripke system ${\mathfrak K}$ based on the set of worlds $K$ as before.
Set $\psi: \A\to \mathfrak{F}_{\mathfrak K}$ by
$\psi_1(p)=(f_l)$ such that if $l=(\Delta, \Gamma)\in K$ is  saturated in $\A_k$,
and $M_k=dim \A_k$, then for $x\in V_k=\bigcup_{p\in G_n}{}^{\alpha}M_k^{(p)}$,
$$f_k(x)=1\Longleftrightarrow {\sf s}_{x\cup (Id_{M_k\sim \alpha)}}^{\A_n}p\in \Delta.$$
Let  $k_0=(\Delta_{\omega}, \Gamma_{\omega})$ be provided by lemma \ref{Ono}, that is a Henkin complete 
saturated extension of $(\Delta_0, \Gamma_0)$
in $\B_1$, then $\psi,$ $k_0$ and $Id$ are as desired. 

\end{demo}

\begin{remark} 
The proof works verabtim when $G$ is stronly rich.
 Theorem \label{Ono2} is the special case of MAIN proved in part 1, when the first pair is just $(\Theta, \bot)$. This weaker version suffices to prove 
interpolation, see the proof of theorem ? in part 3.
\end{remark}

\section{Presence of diagonal elements}

All results, in Part 1, up to the previous theorem,  are proved in the absence of diagonal elements.
Now lets see how far we can go if we have diagonal elements. 
Considering diagonal elements, as we shall see, turn out to be problematic but not hopeless.

Our representation theorem has to respect diagonal elements, 
and this seems to be an impossible task with the presence of infinitary substitutions, 
unless we make a compromise that is, from our point of view, acceptable.
The interaction of substitutions based on infinitary transformations, 
together with the existence of diagonal elements tends to make matters `blow up'; indeed this even happens in the classical case,
when the class of (ordinary) set algebras ceases to be closed under ultraproducts \cite{S}. 
The natural thing to do is to avoid those infinitary substitutions at the start, while finding the interpolant possibly using such substitutions.
We shall also show that in some cases the interpolant has to use infinitary substitutions, even if the original implication uses only finite transformations.

So for an algebra $\A$, we let $\Rd\A$ denote its reduct when we discard infinitary substitutions. $\Rd\A$ satisfies 
cylindric algebra axioms.

\begin{theorem}\label{main3}
Let $\alpha$ be an infinite set. Let $G$ be a semigroup on $\alpha$ containing at least one infinitary transformation. 
Let $\A\in GPHAE_{\alpha}$ be the free $G$ algebra generated by $X$, and suppose that $X=X_1\cup X_2$.
Let $(\Delta_0, \Gamma_0)$, $(\Theta_0, \Gamma_0^*)$ be two consistent theories in $\Sg^{\Rd\A}X_1$ and $\Sg^{\Rd\A}X_2,$ respectively.
Assume that $\Gamma_0\subseteq \Sg^{\A}(X_1\cap X_2)$ and $\Gamma_0\subseteq \Gamma_0^*$.
Assume, further, that  
$(\Delta_0\cap \Theta_0\cap \Sg^{\A}X_1\cap \Sg^{\A}X_2, \Gamma_0)$ is complete in $\Sg^{\Rd\A}X_1\cap \Sg^{\Rd\A}X_2$. 
Then there exist $\mathfrak{K}=(K,\leq \{X_k\}_{k\in K}\{V_k\}_{k\in K}),$ a homomorphism $\psi:\A\to \mathfrak{F}_K,$
$k_0\in K$, and $x\in V_{k_0}$,  such that for all $p\in \Delta_0\cup \Theta_0$ if $\psi(p)=(f_k)$, then $f_{k_0}(x)=1$
and for all $p\in \Gamma_0^*$ if $\psi(p)=(f_k)$, then $f_{k_0}(x)=0$.
\end{theorem}
\begin{demo}{Proof}
The first half of the proof is almost identical to that of  lemma \ref{main}. We highlight the main steps, 
for the convenience of the reader, except that we only deal with the case
when $G$ is strongly rich.
Assume, as usual, that $\alpha$, $G$, $\A$ and $X_1$, $X_2$, and everything else in the hypothesis are given.
Let $I$ be  a set such that  $\beta=|I\sim \alpha|=max(|A|, |\alpha|).$
Let $(K_n:n\in \omega)$ be a family of pairwise disjoint sets such that $|K_n|=\beta.$
Define a sequence of algebras
$\A=\A_0\subseteq \A_1\subseteq \A_2\subseteq \A_2\ldots \subseteq \A_n\ldots$
such that
$\A_{n+1}$ is a minimal dilation of $\A_n$ and $dim(\A_{n+1})=\dim\A_n\cup K_n$.We denote $dim(\A_n)$ by $I_n$ for $n\geq 1$. 
The proofs of Claims 1 and 2 in the proof of \ref{main} are the same.

Now we prove the theorem when $G$ is a strongly rich semigroup.
Let $$K=\{((\Delta, \Gamma), (T,F)): \exists n\in \omega \text { such that } (\Delta, \Gamma), (T,F)$$
$$\text { is a a matched pair of saturated theories in }
\Sg^{\Rd\A_n}X_1, \Sg^{\Rd\A_n}X_2\}.$$
We have $((\Delta_0, \Gamma_0)$, $(\Theta_0, \Gamma_0^*))$ is a matched pair but the theories are not saturated. But by lemma \ref{t3}
there are $T_1=(\Delta_{\omega}, \Gamma_{\omega})$, 
$T_2=(\Theta_{\omega}, \Gamma_{\omega}^*)$ extending 
$(\Delta_0, \Gamma_0)$, $(\Theta_0, \Gamma_0^*)$, such that $T_1$ and $T_2$ are saturated in $\Sg^{\Rd\A_1}X_1$ and $\Sg^{\Rd\A_1}X_2,$ 
respectively. Let $k_0=((\Delta_{\omega}, \Gamma_{\omega}), (\Theta_{\omega}, \Gamma_{\omega}^*)).$ Then $k_0\in K,$
and   $k_0$ will be the desired world and $x$ will be specified later; in fact $x$ will be the identity map on some specified 
domain.

If $i=((\Delta, \Gamma), (T,F))$ is a matched pair of saturated theories in $\Sg^{\Rd\A_n}X_1$ and $\Sg^{\Rd\A_n}X_2$, let $M_i=dim \A_n$, 
where $n$ is the least such number, so $n$ is unique to $i$.
Let $${\bf K}=(K, \leq, \{M_i\}, \{V_i\})_{i\in \mathfrak{K}},$$
where $V_i=\bigcup_{p\in G_n, p\text { a finitary transformation }}{}^{\alpha}M_i^{(p)}$ 
(here we are considering only substitutions that move only finitely many points), 
and $G_n$ 
is the strongly rich semigroup determining the similarity type of $\A_n$, with $n$ 
the least number such $i$ is a saturated matched pair in $\A_n$, and $\leq $ is defined as follows: 
If $i_1=((\Delta_1, \Gamma_1)), (T_1, F_1))$ and $i_2=((\Delta_2, \Gamma_2), (T_2,F_2))$ are in $\bold K$, then set
$$i_1\leq i_2\Longleftrightarrow  M_{i_1}\subseteq M_{i_2}, \Delta_1\subseteq \Delta_2, T_1\subseteq T_2.$$ 
We are not yet there, to preserve diagonal elements we have to factor out $\bold K$ 
by an infinite family equivalence relations, each defined on the dimension of $\A_n$, for some $n$, which will actually turn out to be 
a congruence in an exact sense. 
As usual, using freeness of $\A$, we will  define two maps on $\A_1=\Sg^{\Rd\A}X_1$ and $\A_2=\Sg^{\Rd\A}X_2$, respectively;
then those will be pasted 
to give the required single homomorphism.

Let $i=((\Delta, \Gamma), (T,F))$ be a matched pair of saturated theories in $\Sg^{\Rd\A_n}X_1$ and $\Sg^{\Rd\A_n}X_2$, let $M_i=dim \A_n$, 
where $n$ is the least such number, so $n$ is unique to $i$.
For $k,l\in dim\A_n=I_n$, set $k\sim_i l$ iff ${\sf d}_{kl}^{\A_n}\in \Delta\cup T$. This is well defined since $\Delta\cup T\subseteq \A_n$.
We omit the superscript $\A_n$.
These are infinitely many relations, one for each $i$, defined on $I_n$, with $n$ depending uniquely on $i$, 
we denote them uniformly by $\sim$ to 
avoid complicated unnecessary notation.
We hope that no confusion is likely to ensue. We claim that $\sim$ is an equivalence relation on $I_n.$
Indeed,  $\sim$ is reflexive because ${\sf d}_{ii}=1$ and symmetric 
because ${\sf d}_{ij}={\sf d}_{ji};$
finally $E$ is transitive because for  $k,l,u<\alpha$, with $l\notin \{k,u\}$, 
we have 
$${\sf d}_{kl}\cdot {\sf d}_{lu}\leq {\sf c}_l({\sf d}_{kl}\cdot {\sf d}_{lu})={\sf d}_{ku},$$
and we can assume that $T\cup \Delta$ is closed upwards.
For $\sigma,\tau \in V_k,$ define $\sigma\sim \tau$ iff $\sigma(i)\sim \tau(i)$ for all $i\in \alpha$.
Then clearly $\sigma$ is an equivalence relation on $V_k$. 

Let $W_k=V_k/\sim$, and $\mathfrak{K}=(K, \leq, M_k, W_k)_{k\in K}$, with $\leq$ defined on $K$ as above.
We write $h=[x]$ for $x\in V_k$ if $x(i)/\sim =h(i)$ for all $i\in \alpha$; of course $X$ may not be unique, but this will not matter.
Let $\F_{\mathfrak K}$ be the set algebra based on the new Kripke system ${\mathfrak K}$ obtained by factoring out $\bold K$.

Set $\psi_1: \Sg^{\Rd\A}X_1\to \mathfrak{F}_{\mathfrak K}$ by
$\psi_1(p)=(f_k)$ such that if $k=((\Delta, \Gamma), (T,F))\in K$ 
is a matched pair of saturated theories in $\Sg^{\Rd\A_n}X_1$ and $\Sg^{\Rd\A_n}X_2$,
and $M_k=dim \A_n$, with $n$ unique to $k$, then for $x\in W_k$
$$f_k([x])=1\Longleftrightarrow {\sf s}_{x\cup (Id_{M_k\sim \alpha)}}^{\A_n}p\in \Delta\cup T,$$
with $x\in V_k$ and $[x]\in W_k$ is define as above. 

To avoid cumbersome notation, we 
write ${\sf s}_{x}^{\A_n}p$, or even simply ${\sf s}_xp,$ for 
${\sf s}_{x\cup (Id_{M_k\sim \alpha)}}^{\A_n}p$.  No ambiguity should arise because the dimension $n$ will be clear from context.

We need to check that $\psi_1$ is well defined. 
It suffices to show that if $\sigma, \tau\in V_k$ if $\sigma \sim \tau$ and $p\in \A_n$,  
with $n$ unique to $k$, 
then $${\sf s}_{\tau}p\in \Delta\cup T\text { iff } {\sf s}_{\sigma}p\in \Delta\cup T.$$

This can be proved by induction on the cardinality of 
$J=\{i\in I_n: \sigma i\neq \tau i\}$, which is finite since we are only taking finite substitutions.
If $J$ is empty, the result is obvious. 
Otherwise assume that $k\in J$. We recall the following piece of notation.
For $\eta\in V_k$ and $k,l<\alpha$, write  
$\eta(k\mapsto l)$ for the $\eta'\in V$ that is the same as $\eta$ except
that $\eta'(k)=l.$ 
Now take any 
$$\lambda\in \{\eta\in I_n: \sigma^{-1}\{\eta\}= \tau^{-1}\{\eta\}=\{\eta\}\}\smallsetminus \Delta x.$$
This $\lambda$ exists, because $\sigma$ and $\tau$ are finite transformations and $\A_n$ is a dilation with enough spare dimensions.
We have by cylindric axioms (a)
$${\sf s}_{\sigma}x={\sf s}_{\sigma k}^{\lambda}{\sf s}_{\sigma (k\mapsto \lambda)}p.$$
We also have (b)
$${\sf s}_{\tau k}^{\lambda}({\sf d}_{\lambda, \sigma k}\land {\sf s}_{\sigma} p)
={\sf d}_{\tau k, \sigma k} {\sf s}_{\sigma} p,$$
and (c)
$${\sf s}_{\tau k}^{\lambda}({\sf d}_{\lambda, \sigma k}\land {\sf s}_{\sigma(k\mapsto \lambda)}p)$$
$$= {\sf d}_{\tau k,  \sigma k}\land {\sf s}_{\sigma(k\mapsto \tau k)}p.$$
and (d)
$${\sf d}_{\lambda, \sigma k}\land {\sf s}_{\sigma k}^{\lambda}{\sf s}_{{\sigma}(k\mapsto \lambda)}p=
{\sf d}_{\lambda, \sigma k}\land {\sf s}_{{\sigma}(k\mapsto \lambda)}p$$
Then by (b), (a), (d) and (c), we get,
$${\sf d}_{\tau k, \sigma k}\land {\sf s}_{\sigma} p= 
{\sf s}_{\tau k}^{\lambda}({\sf d}_{\lambda,\sigma k}\cdot {\sf s}_{\sigma}p)$$
$$={\sf s}_{\tau k}^{\lambda}({\sf d}_{\lambda, \sigma k}\land {\sf s}_{\sigma k}^{\lambda}
{\sf s}_{{\sigma}(k\mapsto \lambda)}p)$$
$$={\sf s}_{\tau k}^{\lambda}({\sf d}_{\lambda, \sigma k}\land {\sf s}_{{\sigma}(k\mapsto \lambda)}p)$$
$$= {\sf d}_{\tau k,  \sigma k}\land {\sf s}_{\sigma(k\mapsto \tau k)}p.$$
The conclusion follows from the induction hypothesis.
Now $\psi_1$ respects all quasipolyadic equality operations, that is finite substitutions (with the proof as before; 
recall that we only have finite substitutions since we are considering 
$\Sg^{\Rd\A}X_1$)  except possibly for diagonal elements. 
We check those:

Recall that for a concrete Kripke frame $\F_{\bold W}$ based on ${\bold W}=(W,\leq ,V_k, W_k),$ we have
the concrete diagonal element ${\sf d}_{ij}$ is given by the tuple $(g_k: k\in K)$ such that for $y\in V_k$, $g_k(y)=1$ iff $y(i)=y(j)$.

Now for the abstract diagonal element in $\A$, we have $\psi_1({\sf d}_{ij})=(f_k:k\in K)$, such that if $k=((\Delta, \Gamma), (T,F))$ 
is a matched pair of saturated theories in $\Sg^{\Rd\A_n}X_1$, $\Sg^{\Rd\A_n}X_2$, with $n$ unique to $i$, 
we have $f_k([x])=1$ iff ${\sf s}_{x}{\sf d}_{ij}\in \Delta \cup T$ (this is well defined $\Delta\cup T\subseteq \A_n).$
 
But the latter is equivalent to ${\sf d}_{x(i), x(j)}\in \Delta\cup T$, which in turn is equivalent to $x(i)\sim x(j)$, that is 
$[x](i)=[x](j),$ and so  $(f_k)\in {\sf d}_{ij}^{\F_{\mathfrak K}}$.  
The reverse implication is the same.

We can safely assume that $X_1\cup X_2=X$ generates $\A$.
Let $\psi=\psi_1\cup \psi_2\upharpoonright X$. Then $\psi$ is a function since, by definition, $\psi_1$ and $\psi_2$ 
agree on $X_1\cap X_2$. Now by freeness $\psi$ extends to a homomorphism, 
which we denote also by $\psi$ from $\A$ into $\F_{\mathfrak K}$.
And we are done, as usual, by $\psi$, $k_0$ and $Id\in V_{k_0}$.
\end{demo}

Theorem \ref{main2}, generalizes as is, to the expanded structures by diagonal elements. That is to say, we have:

\begin{theorem}\label{main3}
Let $G$ be the semigroup of finite transformations on an infinite set 
$\alpha$ and let $\delta$ be a cardinal $>0$. Let $\rho\in {}^{\delta}\wp(\alpha)$ be such that
$\alpha\sim \rho(i)$ is infinite for every 
$i\in \delta$. Let $\A$ be the free  $G$ algebra with equality generated by $X$ restristed by $\rho$;
 that is $\A=\Fr_{\delta}^{\rho}GPHAE_{\alpha},$
and suppose that $X=X_1\cup X_2$.
Let $(\Delta_0, \Gamma_0)$, $(\Theta_0, \Gamma_0^*)$ be two consistent theories in $\Sg^{\A}X_1$ and $\Sg^{\A}X_2,$ respectively.
Assume that $\Gamma_0\subseteq \Sg^{\A}(X_1\cap X_2)$ and $\Gamma_0\subseteq \Gamma_0^*$.
Assume, further, that  
$(\Delta_0\cap \Theta_0\cap \Sg^{\A}X_1\cap \Sg^{\A}X_2, \Gamma_0)$ is complete in $\Sg^{\A}X_1\cap \Sg^{\A}X_2$. 
Then there exist a Kripke system $\mathfrak{K}=(K,\leq \{X_k\}_{k\in K}\{V_k\}_{k\in K}),$ a homomorphism $\psi:\A\to \mathfrak{F}_K,$
$k_0\in K$, and $x\in V_{k_0}$,  such that for all $p\in \Delta_0\cup \Theta_0$ if $\psi(p)=(f_k)$, then $f_{k_0}(x)=1$
and for all $p\in \Gamma_0^*$ if $\psi(p)=(f_k)$, then $f_{k_0}(x)=0$.
\end{theorem}
\begin{demo}{Proof} $\Rd\A$ is just $\A$.
\end{demo}  

\begin{remark} Theorem \label{Ono2} extends to the presence of diagonal elements by using a simpler version 
of the congruence relation defined above.
For $n\in \omega$, and a pair $i=(\Delta, \Gamma)$ in  $\A_n,$ 
for $k,l\in dim\A_n$, set $k\sim_i l$ iff ${\sf d}_{kl}^{\A_n}\in \Delta.$
\end{remark}

\subsection{Neat Embedings and the finitizability problem for intuitionistic predicate logic}

Let $\alpha$ be an infinite ordinal.  We set:
$$RGA_{\alpha}={\bf SP}\{\F_{\bold K}: \text { $\bold K$ a Kripke system of dimension $\alpha$}\}.$$
We refer to $RGA_{\alpha}$ as the class of representable $G$ algebras. The next theorem is 
analogous to the celebrated so-calle neat embedding theorem of Henkin in 
cylindric algebras. 
For $\alpha<\beta$, let ${\bf S}\Nr_{\alpha}G_{\beta}PHA_{\beta}=\{\A\in G_{\alpha}PHA_{\alpha}: \exists \B\in G_{\beta}
PHA_{\beta}: \A\subseteq \Nr_{\alpha}\B\}.$

Now using the previous lemmas, we prove a neat embedding theorem:

\begin{theorem}\label{nep}
\begin{enumarab}
\item When $G$ is strongly rich (and $\alpha$ is countable) or $G$ consists of all transformations, then $RGA_{\alpha}=GPHA_{\alpha}$
 \item When $G$ consists only of finite transformtations, then 
we have $RGA_{\alpha}=S\Nr_{\alpha}G_{\alpha+\omega}PHA_{\alpha+\omega}.$ 
In particular,  $RGA_{\alpha}$ is a variety.
\end{enumarab}
\end{theorem}
\begin{demo}{Proof} First it is easy, but tedious, to verify soundness, namely,  that $RGA_{\alpha}\subseteq GPHA_{\alpha}.$ 
One just has to check that equational
axioms for $G$ algebras hold in set algebras based on Kripke systems.
The other inclusion of (i) follows from \ref{rep}, since for a class $K\subseteq GPHA_{\alpha}$, we have $\A\in {\bf SP}K$, if and only if,
for all non-zero $a\in \A$, there exists $\B\in K$ and a homomorphism $f:\A\to \B$ such that $f(a)\neq 0$.

We now show (ii).
It suffices to show that for infinite ordinals 
$\alpha<\beta$,
if $\bold K$ is a Kripke system of dimension $\alpha$ and $\bold M$ is  one of dimension $\beta$, then 
${\mathfrak F}_{\bold K}\subseteq \Nr_{\alpha}{\mathfrak F}_{\bold M}.$
Assume that  ${\mathfrak F}_{\bold K}=\{f=(f_k:k\in K); f_k: {^\alpha}X^{(Id)} \to \mathfrak{O}, k\leq k'\implies f_k\leq f_k'\}$
and that ${\mathfrak F}_{\bold M}=\{g=(g_k:k\in K); g_k: {}^{\beta}X^{(Id)} \to \mathfrak{O}, k\leq k'\implies g_k\leq g_k'\}$
are Kripke set algebras based on the given Kripke systems of dimension $\alpha$ and $\beta$, respectively, with operations
defined as usual. Define $\Psi:{\mathfrak F}_{\bold K}\to {\mathfrak F}_{\bold M}$ 
by  $(\psi f)_k=h_k$ where $h_k\upharpoonright\alpha=f_k$ and 
otherwise is equal to the zero element. Then $\Psi$ neatly embeds $\F_{\bold K}$ into $\F_{\bold M}$, that is it 
embeds $\F_{\bold K}$ into $\Nr_{\alpha}\F_{\bold M}$.
This shows that $RGA_{\alpha}\subseteq S\Nr_{\alpha}G_{\alpha+\omega}PHA_{\alpha+\omega}$. 

The converse inclusion (which is basically an algebraic version of a completeness theorem)
is slighly more difficult.  Assume that 
$\A\subseteq \Nr_{\alpha}\B$. Let $\B'=\Sg^{\B}A$.
Then $\B'\in G_{\alpha+\omega}PHA_{\alpha+\omega}$. This follows from the following reasoning. Let $\beta=\alpha+\omega$.
Then $|\beta\sim \alpha|\geq \omega$, and $\Delta a\subseteq \alpha$ for all $a\in A$, it follows by a simple inductive argument that 
for all $y\in \Sg^{\B}A$, $|\beta\sim \Delta y|\geq \omega$. 
Hence by \ref{rep}, $\B'$ is representable, and so is $\A$, for $\Nr_{\alpha}\B'$ is representable and a subalgebra of a representable algebra is 
representable. To show that $\Nr_{\alpha}\B'$ is representable, it suffices to take
${\mathfrak F}_{\bold K}=\{f=(f_k:k\in K); f_k: {^\beta}X_k^{(Id)} \to \mathfrak{O}, k\leq k'\implies f_k\leq f_k'\}$
and show that $\Nr_{\alpha}\F_{\bold K}$ is representable. But clearly the latter is isomorphic to 
$\F^{\alpha}_{\bold K}=\{(f_k\upharpoonright \alpha:k\in K): f_k\upharpoonright \alpha:{}^{\alpha}X_k^{(Id)}\to \mathfrak{O}, k\leq k', f_k\leq f_k'\}.$

Now we prove something stronger than what is required at the end of item (ii). We 
show that for infinite ordinals $\alpha<\beta$, ${\bf S}\Nr_{\alpha}G_{\beta}PHA_{\beta}$ is a variety. 
It is closed under forming subalgebras  by definition. It is closed under products, since 
for any system of algebras $(\B_i: i\in I)$ we have $\prod_{i\in I} \Nr_{\alpha}\B_i=\Nr_{\alpha}\prod_{i\in I}\B_i$, via the map
$(a_i:i\in I)\mapsto (a_i:i\in I)$ which is evidently one to one, surjective and a homomorphism.

We show that the $S\Nr_{\alpha}G_{\beta}PHA_{\beta}$  is also closed under homomorphic images. 
Let $\A\subseteq \Nr_{\alpha}\B$, and $I$ an ideal of $\A$. We will show that that homomorphic image $\A/I$ of $\A$
is in ${\bf S}\Nr_{\alpha}G_{\beta}PHA_{\beta}$. 
Let $J$ be the ideal generated by $I$ in $\B$. Then $J\cap A=I$.
Define $f:\A/I\to \Nr_{\alpha}(\A/J)$ by $f(x/I)=x/J$, 
then $f$ is an embedding, and we are done. 

\end{demo}

\subsection{The Finitizability for Heyting Polyadic algebras with and without equality}

Item (1) in theorem \ref{nep} says that there is a finite schema axiomatizable variety consisting of subdirect products 
of set algebras based on Kripke systems (models) of a natural extension of ordinary predicate intuitionistic logic. 
This schema can be strictly cut to be a finite axiomatization, that can be easily implemented, using the methods of Sain 
in  \cite{S} where an analagous 
result addressing first order (classical) logic without equality is proved.

The following theorem is the algebraic version of weak completeness see theorem ? in part 3. 
The theorem for classical algebras and its logical consequences are proved in 
\cite{S}. Almost the same proof works in our intuitionistic context (note that (i) and (ii)) are obvious. 
We omit the proof, which can be easily recovered from \cite{S} using fairly standard machinery of algebraic logic dealing 
with the ``bridge" between universal algebraic results and their logical 
counterparts.

In the presence of diagonal elements we have:

\begin{theorem} Let $G$ be rich, then 
we have \begin{enumroman}
\item $GPHAE_{\omega}$ consists of subreducts of full polyadic equality 
Heyting algebras.
\item ${\bf S}\Rd_{ca}GPHAE_{\omega}=RGA_{\omega}.$ 
In other words, the class of representable algebras of dimension $\omega$
coincides with the class of subreducts of $GPHAE_{\omega}$

\item ${\bold H}GPHAE_{\omega}$ is a finitely based variety. 
\end{enumroman}
\end{theorem}
\begin{demo}{Proof} \cite{S}.
\end{demo}
\begin{theorem} For rich $G$, $GPHAE_{\alpha}$ is not closed under ultraproducts
\end{theorem}
In the absence of diagonal elements, a stronger form of the above theorem holds, for in this case $GPHE_{\omega}$ is itself is a finitely 
based variety. In fact in thi case we have:

\begin{theorem} Let $G$ be strongly rich, then  we have \begin{enumroman}
\item $GPHA_{\omega}$ consists of subreducts of full polyadic 
Heyting algebras.
\item ${\bf S}\Rd_{ca}GPHAE_{\omega}=RGH_{\omega}.$ 
In other words, the class of representable algebras of dimension $\omega$
coincides with the class of subreducts of $GPHAE_{\omega}$

\item $RGH_{\omega}$ is a finitely based variety that has $SUPAP$. 
\end{enumroman}
\end{theorem}
\begin{demo}{Proof} \cite{S}.
\end{demo}

In a forthcoming publication  for different semigroups $G$ , we introduce an intuitionistic  
logic (with equality) ${\mathfrak L}_G ({\mathfrak L}_G^{=})$, 
where formulas have infinite length extending predicate intuitionistic logic. 

We show using our algebraic results
that several such logics are complete and have the interpolation properties, while others do not, pending on $G$.
In the presence of equality, our positive results are weaker.

When we deal with those logics whoe atomic formulas exhaust the set of all variables; no variables lie outside formulas,
we discover that the borderline, that turns negative results to positive ones are 
the presence of infinitary transformations, at least two, where one is injective and not onto and the other is its left inverse. 
These transformations  move infinitely 
many variables, and permits one to form dilations in spare dimensions, an essential feature in forming succesive extensions of theories.
We show that the presence of such transfornations is also necessary, in the sense that their presence cannot be dispensed with. If only finite transformations are available, 
then, without imposing any conditions on dimension sets of elements,  representability and amalgamation
fail in a strong sense.

\end{document}